\documentclass{amsart}

\usepackage{amsmath,amsthm,amssymb,amscd,mathtools}
\usepackage{amsfonts}
\usepackage{rotating}
\usepackage{euscript}
\usepackage{pst-node}
\usepackage{epsfig,verbatim}
\def\be{\begin{equation}}
\def\ee{\end{equation}}

\def\C{{\mathbb C}}

\def\P{{\mathbb P}}

\def\phi{{\varphi}}
\def\v{{\varepsilon}} 

\def\deg{{\rm deg\,}}

\def\bp{\begin{proposition}}
\def\ep{\end{proposition}}

\def\bt{\begin{theorem}}
\def\et{\end{theorem}}
\def\br{\begin{remark}}
\def\er{\end{remark}}
\def\be{\begin{equation}}
\def\bee{\begin{equation*}}
\def\l{\label}

\def\m{\mu}
\def\ee{\end{equation}}
\def\eee{\end{equation*}}
\def\bl{\begin{lemma}}
\def\el{\end{lemma}}
\def\bc{\begin{corollary}}
\def\ec{\end{corollary}}
\def\pr{\noindent{\it Proof. }}

\def\bd{\begin{definition}}
\def\ed{\end{definition}}

\mathtoolsset{showonlyrefs}

\newtheorem{theorem}{Theorem}[section]
\newtheorem{lemma}{Lemma}[section]
\newtheorem{definition}{Definition}[section]
\newtheorem{corollary}{Corollary}[section]
\newtheorem{proposition}{Proposition}[section]
\newtheorem{remark}{Remark}[section]

\begin{document}
\title[Functions sharing the measure of maximal entropy]{
On rational functions sharing the measure of maximal entropy  
}
 
\author{F. Pakovich}
\thanks{
This research was  supported by ISF Grant No. 1432/18}
\address{Department of Mathematics, Ben Gurion University of the Negev, Israel}
\email{
pakovich@math.bgu.ac.il}

\begin{abstract} 
We show that describing rational functions $f_1,$ $f_2,$ $\dots,f_n$ sharing the 
measure of maximal entropy reduces to describing solutions of the functional equation 
$A\circ X_1=A\circ X_2=\dots=A\circ X_n$ in rational functions. 
We also provide some results about solutions of this equation.

\end{abstract}

\maketitle

\begin{section}{Introduction}

 Let $f$ be a rational function of degree $d\geq 2$ on $\C\P^1$. It was proved by 
 Friere, Lopes, and Ma\~{n}\'e (\cite{flm}), and independently by  Lyubich (\cite{l}) 
that there exists a unique probability measure $\mu_f$ on $\C\P^1$, which is invariant under $f$, has support equal to the Julia set $J(f)$ of $f$, and achieves maximal entropy 
$\log d$ among all $f$-invariant probability measures.
In this note, we study rational functions 
sharing the measure of maximal entropy, that is rational functions 
$f$ and $g$ such that $\mu_f=\mu_g,$ and more generally rational functions  $f_1,$ $f_2,$ $\dots,f_n$
such that $\mu_{f_1}=\mu_{f_2}=\dots =\mu_{f_n}.$
We assume 
that considered functions 
are {\it non-special} in the following sense: they are  
neither Latt\`es maps nor conjugate to $z^{\pm n}$ or $\pm T_n.$

In case if $f$ and $g$ are polynomials, the  condition $\mu_f=\mu_g$ is equivalent to the condition 
$J(f)=J(g).$ 
In turn, for non-special polynomials $f$ and $g$ the equality $J(f)=J(g)=J$ holds if and only if there exists 
a polynomial $h$ such that $J(h)=J$ and 
\be \l{for} f=\eta_1\circ h^{\circ s}, \ \ \  g=\eta_2\circ h^{\circ t}\ee for some integers $s,t\geq 1$ and 
rotational symmetries $\eta_1,$ $\eta_2$ of $J$  (see \cite{a2}, \cite{sh}, and also 
 \cite{a1},  \cite{b}-\cite{be2} for other related results).
Note that a similar conclusion remains true if instead of the condition $J(f)=J(g)$ one were to assume only that $f$ and $g$ share a completely invariant compact set in $\C$ (see \cite{p1}). 
Note also  that in the polynomial case any of the conditions  $J(f)=J(g)$ and 
\eqref{for} is equivalent to the condition that 
\be \l{cond} f^{\circ k}=\eta \circ g^{\circ l},\ee for some integers $k,l\geq 1$ and
 M\"obius transformation $\eta$  such that $\eta(J(g))=J(g).$

Since $\mu_f=\mu_{f^{\circ k}},$ the equality 
 $\mu_f=\mu_g$ holds whenever $f$ and $g$ share an iterate, that is satisfy \be \l{ite}  f^{\circ k}=g^{\circ l}\ee for some integers  $k,l\geq 1.$ Moreover,  $\mu_f=\mu_g$ whenever $f$ and $g$ commute. 
However, the latter condition in fact is a particular case of the former one, since 
 non-special commuting $f$ and $g$ always satisfy \eqref{ite} by the result of Ritt (\cite{r2}).
Note that in distinction with the polynomial case rational solutions of \eqref{ite} not necessarily have the form \eqref{for} (see \cite{r2}, \cite{pp}).

The problem of describing rational functions $f$ and $g$ with $\mu_f=\mu_g$  can be expressed in 
algebraic terms. 
Specifically,  the results of Levin (\cite{lev}) and Levin and Przytycki
(\cite{lp}) imply that for non-special  
$f$ and $g$ the equality $\mu_f=\mu_g$ holds if and only if 
some of their iterates 
$F=f^{\circ k}$ and $G=g^{\circ l}$ satisfy the system of functional equations 
\be \l{eq} F\circ F=F\circ G, \ \ \ \ G\circ G =G\circ F\ee 
(see \cite{ye} for more detail).

Examples of rational functions $f$, $g$ with $\mu_f=\mu_g$, which do not have the form \eqref{cond}, 
were constructed by Ye in the paper \cite{ye}.
These examples are based on the following remarkable observation: if $X$, $Y$, and $A$ are  rational functions such that 
\be \l{fe} A\circ X=A\circ Y,\ee then the functions \be \l{sys} F=X\circ A, \ \ \  \ G=Y\circ A\ee
satisfy \eqref{eq}.    
The simplest examples of solutions of \eqref{fe} can be obtained from rational functions satisfying $A\circ \eta=A$
for some M\"obius transformation $\eta$, 
by setting  
\be \l{mu} X=\eta\circ Y.\ee In this case, the corresponding solutions of \eqref{eq} have the form \eqref{cond}. However, other solutions of \eqref{fe} also exist, allowing to construct solutions of \eqref{eq} which do not have the form \eqref{cond}. 

Roughly speaking, the main result of this note states that in fact {\it all} solutions of \eqref{eq} 
can be obtained from solutions of \eqref{fe}. More generally, the following statement holds.

\bt \l{t1} Let $f_1,f_2,\dots, f_n$  be non-special rational functions of degree at least two on $\C\P^1$.   Then they share the measure of maximal entropy if and only if 
some of their iterates $F_1,F_2,\dots,F_n$
can be represented in the form 
\be \l{func} F_1=X_1\circ A, \ \ \  \ F_2=X_2\circ A, \ \ \ \ \dots \ \ \ \ ,
 F_n=X_n\circ A,\ee 
where $A$ and $X_1,X_2,\dots,X_n$  are  rational functions
such that 
\be \l{sys2} A\circ X_1=A\circ X_2=\dots =A\circ X_n \ee
and
$\C(X_1,X_2,\dots,X_n)=\C(z).$ 
\et 

Theorem \ref{t1} shows that ``up to iterates" 
describing pairs of rational functions $f$ and $g$ 
with $\mu_f=\m_g$ 
reduces to describing solutions of \eqref{fe}. In particular, since polynomial solutions of \eqref{fe} satisfy \eqref{mu}, 
we immediately recover the result that  polynomials $f,g$ with $\mu_f=\m_g$ satisfy \eqref{cond}. 
Nevertheless, the problem of describing solutions of \eqref{fe} for arbitrary rational $A$, $X$, $Y$ 
is still widely open. In fact, a complete description    of solutions of \eqref{fe} 
is obtained only in the case where $A$ is a polynomial (while $X$ and $Y$ can be arbitrary rational functions)
in the paper \cite{az} by Avanzi and Zannier. The  approach of \cite{az} is based on describing polynomials  $A$ for which the genus of an irreducible algebraic curve 
\be \l{cucu} C_A: \ \frac{A(x)-A(y)}{x-y}=0\ee
 is  zero, and analyzing  situations where $C_A$ is reducible but has a component of genus zero. Although the same strategy can be applied 
to an arbitrary rational function $A$, both its stages become much more complicated and no general results are known to date.

Note  that the problem of describing solutions of equation \eqref{fe} for rational $A$ and 
  {\it meromorphic}  on the complex plane $X$, $Y$ was posed in the paper of Lyubich and Minsky  (see \cite{lm}, p. 83) 
in the context of studying the action of rational functions on the ``universal space'' of non-constant functions meromorphic on $\C$. In algebraic terms, the last problem is equivalent to describing  
rational functions $A$ such that \eqref{cucu} has a component of genus zero {\it or one}.

 Theorem \ref{t1} implies an interesting corollary, concerning dynamical characteristics of rational functions sharing the measure of maximal entropy. 
Recall that the {\it multiplier spectrum} of a rational function $f$ of degree $d$ is a function which assigns to each $s\geq 1$ the unordered list of multipliers at all $d^s+1$ fixed points of $f^{\circ s}$ taken with appropriate multiplicity.  
Two rational functions  are called {\it isospectral} if they have the same   multiplier spectrum.

\bc \l{c1} If non-special rational functions $f_1,f_2,\dots, f_n$  of degree at least two share the measure of maximal entropy, then 
some of their iterates $F_1,F_2,\dots,F_n$  are isospectral.
\ec

The rest of this note is organized as follows. In the second section, we prove Theorem \ref{t1} and Corollary \ref{c1}. Then, in the third section, we prove two results concerning equation 
 \eqref{fe} and system \eqref{sys2}. The first result states that if the curve 
$C_A$ is irreducible and rational functions $X$, $Y$ provide a generically one-to-one parametrization of $C_A$, 
then  $X=Y\circ \eta$ for some involution $\eta\in Aut(\C\P^1).$
The second result states that if $A$ and $X_1,X_2,\dots, X_n$ are rational functions such that 
\eqref{sys2} holds and $X_1,X_2,\dots, X_n$ are distinct, then 
$n\leq \deg A$, and $n= \deg A$ only if  
the Galois closure of the field extension  $\C(z)/\C(A)$ has genus zero or one.
In fact, we prove these results in the more general setting, allowing the functions 
 $X$, $Y$ and $X_1,X_2,\dots,X_n$ to be meromorphic on $\C.$

\end{section}

\begin{section}{Functions sharing the measure of maximal entropy}
In this section, we deduce Theorem \ref{t1} and Corollary \ref{c1} from the 
criterion \eqref{eq} and 
the following four lemmas.

\bl \l{02}
Let $A_1,A_2,\dots, A_n$ and $Y_1,Y_2,\dots, Y_n$ be  rational functions such that
\be \l{syska} A_i\circ Y_1=A_i\circ  Y_2=\dots =A_i\circ Y_n, \ \ \ i=1,\dots n,\ee and \be \l{es} \C(A_1,A_2,\dots,A_n)=\C(z).\ee Then 
\be \l{rava} Y_1=Y_2=\dots =Y_n.\ee
\el 
\pr  
By \eqref{es}, there exists a rational function $P\in \C(z_1,z_2,\dots,z_n)$ such that 
\be \l{r} z=P(A_1,A_2,\dots,A_n),\ee
implying that 
\be \l{xre} Y_j=P(A_1\circ Y_j,A_2\circ Y_j,\dots,A_n\circ Y_j), \ \ \ \ 1\leq j \leq n.\ee Now \eqref{rava} follows from \eqref{xre} and \eqref{syska}. \qed

\bl \l{koro} 
Let $F_1,F_2,\dots, F_n$ be rational functions such that
\be \l{s1} F_i\circ F_1=F_i\circ  F_2=\dots =F_i\circ F_n, \ \ \ i=1,\dots n.\ee Then 
 there exist rational functions $A$ and  $X_1,X_2,\dots, X_n$   such that 
\be \l{se} F_i=X_i\circ A, \ \ \ i=1,\dots n,\ee
\be \l{sd} \C(X_1,X_2,\dots,X_n)=\C(z),\ee and 
\be \l{s2} A\circ X_1=A\circ X_2=\dots =A\circ X_n.\ee
\el
\pr By the L\"uroth theorem,   
$$\C(F_1,F_2,\dots,F_n)=\C(A)$$ for some rational function $A$, implying that equalities 
 \eqref{se}  hold for some rational functions $X_1,X_2,\dots, X_n$ satisfying \eqref{sd}.
Substituting now \eqref{se} in \eqref{s1}  
we see that 
$$X_i\circ (A\circ X_1)=X_i\circ (A\circ X_2) =\dots =X_i\circ (A\circ X_n), \ \ \ i=1,\dots n.$$ 
Applying now Lemma \ref{02} to the last system 
we obtain \eqref{s2}. \qed

\bl \l{01} Let $A$ and $B$ be rational functions such that the equality $$A\circ A=A\circ B$$ holds. Then 
$$A^{\circ l}\circ A^{\circ l}=A^{\circ l}\circ B^{\circ l}$$
for any 
$l\geq 1$.
\el \pr The proof is by induction on $l$. Assuming that the lemma is true for $l=k$, we have: 
$$A^{\circ (k+1)}\circ B^{\circ (k+1)}=A^{\circ k}\circ (A\circ B)\circ B^{\circ k}=A^{\circ k}\circ A^{\circ 2}\circ B^{\circ k}=$$
$$=A^{\circ 2}\circ A^{\circ k}\circ B^{\circ k}=A^{\circ 2}\circ A^{\circ 2k}=A^{2k+2}.\eqno{\Box}$$

\bl \l{od} Let $d_i\geq 2$, $1\leq i \leq n,$ and  $n_{i,j}\geq 1,$ $1\leq i,j\leq n,$ $i\neq j,$ 
be integers such that 
$$d_i^{n_{i,j}}=d_j^{n_{j,i}},\ \ \ 1\leq i,j\leq n, \ \ \ i\neq j.$$
Then there exist integers $l_i\geq 1,$ $1\leq i \leq n,$
such that 
$$d_1^{l_1}=d_2^{l_2}=\dots =d_n^{l_n}.$$
\el
\pr The proof is by induction on $n$. For $n=2$, we obviously can set $$l_1=n_{1,2},\ \ \ \ l_2=n_{2,1}.$$
Assuming that the lemma is true for $n=k$, we can find integers $a_i,$ $1\leq i \leq k,$ and $b_i,$ $2\leq i \leq k+1,$
such that 
$$d_1^{a_1}=d_2^{a_2}=\dots =d_{k}^{a_{k}}$$
and 
$$d_2^{b_2}=d_3^{b_3}=\dots =d_{k+1}^{b_{k+1}},$$ implying that 
$$d_1^{a_1b_2}=d_2^{a_2b_2}=\dots =d_k^{a_{k}b_2}$$
and 
$$d_2^{b_2a_2}=d_3^{b_3a_2}=\dots =d_{k+1}^{b_{k+1}a_2}.$$
Therefore, 
$$d_1^{a_1b_2}=d_2^{b_2a_2}=d_3^{b_3a_2}=\dots =d_{k+1}^{b_{k+1}a_2},$$
and hence the lemma is true for $n=k+1.$ \qed 
%we can set  $$l_1=a_1b_2, \ \ \ l_2=b_2a_2, \ \ \ l_3=b_3a_2, \ \ \  \dots \ \ \ ,  l_{k+1}=b_{k+1}a_2. \eqno{\Box}$$
\vskip 0.2cm

\noindent{\it Proof of Theorem \ref{t1}.}
For any rational functions  $A$ and $X_1,X_2,$ $\dots, X_n$ satisfying \eqref{sys2} the corresponding functions \eqref{func} satisfy system \eqref{s1}. In particular, for any pair $i,j$  $1\leq i,j\leq n,$ $i\neq j,$ the equalities 
$$F_i\circ F_i=F_i\circ F_j, \ \ \ \ F_j\circ F_j=F_j\circ F_i, \ \ \ \ 1\leq i,j\leq n,$$ hold, implying that   the functions $f_i$, $f_j$ share the measure of maximal entropy. Therefore,
all $f_1, f_2, \dots, f_n$ share the measure of maximal entropy.

In the other direction, if $\mu_{f_1}=\mu_{f_2}=\dots =\mu_{f_n},$ 
then using the criterion \eqref{eq} we can find  integers $n_{i,j},$ $1\leq i,j\leq n,$ $i\neq j,$ 
such that 
\be \l{ekr} f_i^{\circ n_{i,j}}\circ f_i^{\circ n_{i,j}}=f_i^{\circ n_{i,j}}\circ f_j^{\circ n_{j,i}}, \ \ \ \ \ 
f_j^{\circ n_{j,i}}\circ f_j^{\circ n_{j,i}}=f_j^{\circ n_{j,i}}\circ f_i^{\circ n_{i,j}}.\ee 
Suppose first that \be \l{usa} \deg f_1=\deg f_2=\dots =\deg f_n.\ee 
Then \eqref{usa} and \eqref{ekr} imply that $n_{i,j}=n_{j,i},$ $1\leq i,j\leq n.$
Applying now Lemma \ref{01} to \eqref{ekr}, we see that for any integer number $M$ divisible by all the numbers $n_{i,j},$ $1\leq i,j\leq n,$ the equalities 
$$f_i^{\circ M}\circ f_i^{\circ M}=f_i^{\circ M}\circ f_j^{\circ M}, \ \ \ \ \ 
1\leq i,j\leq n,$$ hold. Thus, the functions $F_i=f_i^{\circ M},$ $ 1\leq i \leq n,$ satisfy system \eqref{s1}, implying by Lemma \ref{koro}  that equalities \eqref{se}, \eqref{sd}, and \eqref{s2} hold.

For arbitrary rational functions $f_1, f_2, \dots, f_n$ sharing the measure of maximal entropy, we still can write system \eqref{ekr}, implying that 
$$(\deg f_i)^{n_{i,j}}=(\deg f_j)^{n_{j,i}}, \ \ \  1\leq i,j\leq n,\ \ \  i\neq j.$$ Applying Lemma \ref{od}, we can find  $l_{i},$ $1\leq i\leq n,$ 
such that the rational functions $f_i^{\circ l_i},$ $1\leq i \leq n,$ 
have  the same degree. Since these functions along with the functions $f_1, f_2, \dots, f_n$ share 
the measure of maximal entropy, we can write system \eqref{ekr} for these functions. Using now the already proved part of the theorem, we conclude that there exist  
 $m_{i},$ $1\leq i\leq n,$ 
such that the rational functions $ F_i=f_i^{\circ m_i},$ $1\leq i \leq n,$ 
satisfy  \eqref{s1}, implying \eqref{se}, \eqref{sd}, and \eqref{s2}. \qed 

\vskip 0.2cm
\noindent{\it Proof of Corollary \ref{c1}.}
The corollary follows from the statement of the theorem and the fact that for any   rational functions $U$ and $V$ the rational functions $U \circ V$ and
$V \circ U$ are isospectral (see \cite{rec}, Lemma 2.1). \qed

\end{section}

\begin{section}{Functional equation $A(\phi)=A(\psi)$.}
Equation \eqref{fe} is a particular case of the functional equation 
\be \l{ab} A\circ X=B\circ Y,\ee 
which, under different assumptions on $A,B$ and $X,Y$, has been studied in many papers (see e. g. \cite{an}, \cite{az2}, \cite{bilu}, \cite{f1}, \cite{nv}, \cite{pak},  \cite{pakamer}, \cite{cur}, \cite{r1}).
Nevertheless, to our best knowledge precisely equation \eqref{fe} was the subject of only  two papers. One of them is the paper of Avanzi and Zannier cited in the introduction. The other one 
is the paper \cite{r4} by Ritt, written eighty years earlier,
where some partial results were obtained.  
In particular, Ritt observed that  solutions of \eqref{fe} with $X\neq Y$ 
can be obtained using  
finite subgroups of $Aut(\C\P^1)$ as follows. Let $\Gamma$ be a finite subgroup of $Aut(\C\P^1)$ and 
$\theta_{\Gamma}$ its invariant function, that is a rational function such that $\theta_{\Gamma}(x)=
\theta_{\Gamma}(y)$ if and only if $y=\sigma(x)$ for some $\sigma\in \Gamma.$ Then for any subgroup $\Gamma'\subset \Gamma$ the equality 
\be \l{eqq} \theta_{\Gamma}=\psi \circ \theta_{\Gamma'}\ee holds for some $\psi\in \C(z)$, implying that   
$$\psi\circ \theta_{\Gamma'}=\psi\circ (\theta_{\Gamma'}\circ \sigma)$$ for every $\sigma\in \Gamma$.  
Nevertheless, $\theta_{\Gamma'}\neq \theta_{\Gamma'}\circ \sigma$ unless $\sigma\in \Gamma'.$
For example, for the dihedral group $D_{2n}$, generated by 
$z\rightarrow 1/z$ and $z\rightarrow \v z,$ where $\v=e^{\frac{2\pi i}{n}}$, and its subgroup $D_2$ 
equality \eqref{eqq} takes the form
$$\frac{1}{2}\left(z^n+\frac{1}{z^n}\right)=T_n\circ \frac{1}{2}\left(z+\frac{1}{z}\right)$$ 
giving rise to the solution
\be \l{xru} T_n\circ \frac{1}{2}\left(z+\frac{1}{z}\right)=T_n\circ \frac{1}{2}\left(\v z+\frac{1}{\v z}\right)\ee 
of \eqref{fe} not satisfying to \eqref{mu}. 
Ritt also constructed  solutions of \eqref{fe} using rational functions arising from the formulas for the period transformations  of the Weierstrass functions $\wp(z)$ for lattices with symmetries of order greater than two.

In this note, we do not make an attempt to obtain an explicit classification of solutions of \eqref{fe} in spirit of \cite{az}. Instead, we prove two general results which emphasize the role of symmetries 
 in the problem.

\bt \l{t2} Let $A$ be a rational function and $\phi,\psi$ distinct functions meromorphic  on $\C$ 
such that $$A\circ \phi=A\circ \psi.$$ Assume in addition that 
the algebraic curve $C_A$ is irreducible. Then the desingularization  $R$ of $C_A$
has genus zero or one and there exist holomorphic functions 
$\phi_1:\, R\rightarrow \C\P^1$, $\psi_1:\, R\rightarrow \C\P^1$ and $h:\,\C \rightarrow R$  
 such that 
\be \l{phi} \phi=\phi_1 \circ h, \ \ \ \ \  \psi=\psi_1  \circ h,\ee 
and   the map from $R$ to $C_A$ given by 
$ z\rightarrow (\phi_1(z), \psi_1(z))$ is generically one-to-one.
Moreover, 
\be\l{sati}\phi_1=\psi_1\circ \eta\ee for some involution $\eta:R\rightarrow R$.
\et
\pr 
The first conclusion of the theorem holds for any  parametrization of an al\-geb\-raic curve by functions meromorphic on $\C$ (see e. g. \cite{bn}, Theorem 1 and Theorem 2), so we only must show the existence of an involution $\mu$ 
satisfying \eqref{sati}. 

Since the equation of $C_A$ is invariant under the exchange of variable,  along with the meromorphic parametrization $z\rightarrow (\phi_1,\psi_1)$ the curve $C_A$  admits the meromorphic parametrization $z\rightarrow (\psi_1,\phi_1)$. Since the desingularization  $R$ is defined up to an automorphism, it follows now from the first part 
of the theorem that 
$$\phi_1=\psi_1\circ \eta, \ \ \ \ \ \ \ \psi_1=\phi_1\circ \eta$$ for some $\eta\in Aut(R)$, implying that 
\be \l{sl} \phi_1=\phi_1\circ (\eta\circ \eta), \ \ \ \ \ \ \ \psi_1=\psi_1\circ (\eta\circ \eta).\ee
Finally, $\eta \circ \eta=z$ since otherwise \eqref{sl} contradicts to the condition that 
the map 
$ z\rightarrow (\phi_1(z), \psi_1(z))$ is generically one-to-one. \qed

\bt \l{t3} Let $A$ be a rational function of degree $d$ and $\phi_1,\phi_2,\dots,\phi_n$ distinct meromorphic functions on $\C$ such that 
\be \l{hold} A\circ \phi_1=A\circ \phi_2=\dots =A\circ \phi_n.\ee Then $n\leq d$. Moreover, if $n=d$, then 
the Galois closure of the field extension  $\C(z)/\C(A)$ has genus zero or one. 
\et
\pr Since for any $z_0\in \C\P^1$  the preimage $A^{-1}(z_0)$ contains at most 
$d$ distinct points, if \eqref{hold} holds for $n>d$, then for every $z\in \C\P^1$ at most $d$ of the values $\phi_1(z),\phi_2(z),\dots,\phi_n(z)$ are distinct, implying that  at most $d$ of the functions  $\phi_1,\phi_2,\dots,\phi_n$ are distinct. 

The second part of the theorem is the ``if'' part of the following criterion  (see \cite{da}, Theorem 2.3). For a rational function $A$ of degree $d$, the Galois closure of the field extension  $\C(z)/\C(A)$ has genus zero or one if and only if there exist $d$ distinct functions $\psi_1, \psi_2, \dots ,\psi_d$ meromorphic on $\C$ such that 
$$ A\circ \psi_1=A\circ \psi_2=\dots =A\circ \psi_d.\eqno{\Box}$$

Note that  rational functions $A$ for which the genus $g_A$ of the Galois closure of the field extension  $\C(z)/\C(A)$ is zero are exactly all possible ``compositional left factors'' of Galois coverings of
$\C\P^1$ by $\C\P^1$ and can be listed explicitly. On the other hand, functions with $g_A=1$
admit a simple geometric description in terms of projections of maps between elliptic curves 
(see \cite{gen}). The simplest examples of rational functions with $g_A\leq 1$ are $z^n,$ $T_n$, $\frac{1}{2}\left(z^n+\frac{1}{z^n}\right)$, and Latt\`es maps.

\end{section}

\bibliographystyle{amsplain}

\end{document}